\documentclass[11pt]{article}
\usepackage[cp1251]{inputenc} 
\usepackage[dvips]{graphicx}
\usepackage[dvips]{epsfig}
\usepackage{amssymb}
\begin{document}

\begin{center}
\textsc{A. R. Mirotin}
\end{center}

\begin{center}
\Large{{\bf  PERTURBATION DETERMINANTS  ON BANACH SPACES AND OPERATOR DIFFERENTIABILITY FOR HIRSCH FUNCTIONAL CALCULUS}}
\end{center}

\vspace{1cm}

\begin{center}
amirotin@yandex.ru
\end{center}

\vspace{1cm}

\textbf{Abstract.}
{\small We consider a perturbation determinant for  pairs of nonpositive (in a sense of Komatsu) operators on  Banach space with nuclear difference  and prove a generalization of the important  formula for the logarithmic derivative of this  determinant.    To this end the  Frechet  differentiability of operator monotonic (negative complete Bernstein) functions of negative and nonpositive operators on Banach spaces is investigated.
The results may be regarded as  a contribution to the Hirsch functional calculus.}

\vspace{1cm}

\textbf{Key wards.}    Perturbation determinant;  nonpositive operator; Hirsch functional calculus; Bernstein function; operator monotonic function;  operator differentiability.
\vspace{1cm}

\vspace{2cm}

{\bf 1. Introduction}

The   perturbation determinant  plays  a very important role in perturbation theory for linear operators.  It was introduced by  M. Krein  in his seminal paper \cite{Krein}   for operators on Hilbert space and is crucial in particular for Krein's theory of trace formulas  and  spectral shift functions. Later this concept received many other applications (see, e.g., \cite[Chapter IV]{GK} where perturbation determinants are channels for the use in operator theory of theorems of the theory of functions). For more recent results on perturbation determinants see, e.g., \cite{CJLS}, \cite{GZ}, \cite{MS}, \cite{MN1}. As was mentioned in \cite[Chapter IV]{GK} a number of relations of the theory of perturbation determinants can
be traced as far back as an old paper by H.~Bateman \cite{Ba} on integral
operators. All aforementioned works deal with Hilbert space operators.

 The  perturbation determinants for the pairs of  generators of strongly continuous semigroups on Banach spaces was considered by the author in \cite{OaMII}.
 In the present article, following \cite{arX2018} we   extend the classical concept of perturbation determinant to a more general setup   of pairs of nonpositive operators  on  Banach space with nuclear difference  and prove a generalization of the important  formula for the logarithmic derivative of this  determinant. The last formula was proved (by another method) in \cite{arX2018} for negative operators   only.  Our main tool is the notion of a Frechet derivative of a Bernstein function of an operator argument.

Yu. Daletski\u{\i} and S.G. Kre\u{\i}n pioneered the study of the problem of differentiability of functions of self-adjoint operators   in  \cite{DK}. Their study  has been motivated
by problems in perturbation theory. Existence of the higher order derivatives  was established
 in \cite{Pel06}.
 Differential calculus for functions of several commuting Hermitian operators in Hilbert spaces was studied in \cite{KPSSII}. For a survey and bibliography of the theory for Hilbert spaces that resulted  see the  article \cite{AP1}.   The case of Banach spaces was considered in  \cite{RM},  \cite{OaM}, and \cite{OaMII}.

 In this paper (bearing in mind the applications to the theory of perturbation determinants on Banach spaces) we investigate operator differentiability of operator monotonic (negative complete Bernstein) functions of negative and nonpositive operators on Banach spaces. Thus we generalize results on operator differentiability  obtained for compactly supported representing measures in \cite{RM} to the case of measures supported in $\mathbb{R}_+.$ These results may also have an independent interest.

The paper may be regarded as  a contribution to a perturbation theory for Hirsch functional calculus.

{\bf 2.  Preliminaries}

In this section we introduce classes of functions and  operators and  briefly  describe  a version of Hirsch functional calculus  we shall use below.

\textbf{Definition 1.} We say that a function $\varphi$ is  \textit{negative complete Bernstein} and write $\varphi\in \mathcal{OM}_-$ if it is holomorphic in $\mathbb{C}\setminus \mathbb{R}_+$, satisfies $\mathrm{Im} w \mathrm{Im}\varphi(w)\geq 0$ for $w\in\mathbb{C}\setminus \mathbb{R}_+,$ and such that the limit $\varphi(-0)$ exists
and is real.

According to \cite[Theorem 6.1]{SSV} this means that $-\varphi(-z)$ is a complete Bernstein function and $\varphi$ has the following integral representation
$$
\varphi(z)=c+bz+\int\limits_{(0,\infty)}\frac{z}{t-z}d\mu(t),\quad z\in \mathbb{C}\setminus (0,+\infty),\eqno(1)
$$
where $c\leq 0, b\geq 0$ and $\mu$ is a unique positive measure such that $\int_{(0,\infty)}d\mu(t)/(1+t)<\infty.$

A lot of  examples of complete Bernstein functions  one can found in \cite{SSV}.

In the sequel unless otherwise stated we  assume for the sake of simplicity  that $c=b=0$ in the integral representation (1) (otherwise one should replace $\varphi(z)$ by $\varphi(z)-c-bz$).

\textbf{Remark 1.} It is known (see, e.g., \cite[Theorem 12.17]{SSV}), that the families of complete Bernstein and positive operator monotone functions coincide. It follows that the families of negative complete Bernstein and negative operator monotone functions also coincide (we say that a real function $\varphi$ on $(-\infty,0]$ is \textit{negative operator monotone} if for every bounded self-adjoint operators $A$ and $B$ on a finite or infinite-dimensional real Hilbert space the inequalities $A\leq B\leq O$ imply $\varphi(A)\leq \varphi(B)$). That is why  we denote
the family of negative complete Bernstein functions by $\mathcal{OM}_-.$

\textbf{Definition 2.} We say that closed, densely defined operator $A$ on a complex Banach space $X$ is \textit{nonpositive (negative)} (in a sense of Komatsu) if
 $(0,\infty)$ is contained in $\rho(A),$ the resolvent set of $A$ (respectively  $[0,\infty)\subset \rho(A)$), and
 $$
 M_A:=\sup_{t>0}\|tR(t,A)\|<\infty
 $$
  (respectively
  $$
 M_A:=\sup_{t>0}\|(1+t)R(t,A)\|<\infty)
 $$
  where $R(t,A) = (tI- A)^{-1}$ stands for the resolvent of an operator $A,$ and $Ix=x$ for all $x\in X.$

We denote by  $\mathcal{NP}(X)$  ($\mathcal{N}(X)$) the class of nonpositive (respectively, negative) operators on the space $X.$

\textbf{Remark 2.} The operator $A$  is nonpositive (negative)  if and only if $-A$ is nonnegative (positive) in a sense of Komatsu \cite{Kom} (see also \cite[Chapter 1]{MS}).  We   deal with  negative  operators instead of positive one because in this form our results are consistent with Bochner-Phillips functional calculus of semigroup generators considered  in \cite{Mir99} -- \cite{MirSF}.
According to the Hille-Yosida Theorem every generator of strongly continuous uniformly bounded semigroup of operators belongs to  $\mathcal{NP}(X).$ It follows that even for the case of Hilbert space the class of nonpositive operators in a sense of Komatsu is wider that the class of self-adjoint    nonpositive operators  in the classical sense.

\textbf{Definition 3}  \cite{BBD}. For any function  $\varphi\in \mathcal{OM}_-$ with representing measure $\mu$ and any $A\in \mathcal{NP}(X)$ we put
$$
 \varphi(A)x=\int\limits_{(0,\infty)}AR(t,A)xd\mu(t) \quad (x\in D(A)) \eqno(2)
$$
(the Bochner integral). This operator is closable (see, e.g., \cite{BBD}) and its closure will be denoted by $\varphi(A)$, too.

It is known \cite{Hirsch72} --\cite{Hirsch76} (see also \cite[Theorem 7.4.6]{MS}) that for $\varphi\in \mathcal{OM}_-$ the operator $\varphi(A)$ belongs to $\mathcal{NP}(X)$ ($\mathcal{N}(X)$) if $A\in \mathcal{NP}(X)$ (respectively, $A\in \mathcal{N}(X)$).

\textbf{Remark 3.} In the Hirsch functional calculus \cite{Hirsch72} --\cite{Hirsch76} (see also \cite{MS}, \cite{Pys}, \cite{TrudyIM}) functions of the form
$$
f(w)=a+\int\limits_{[0,\infty)}\frac{w}{1+ws}d\lambda(s)
$$
 ($a\geq 0,$  $\lambda$ is a unique positive measure such that $\int\limits_{(0,\infty)}d\lambda(s)/(1+s)<\infty$) are applied to nonnegative operators $T$ on Banach spaces via the formula
 $$
f(T)x=ax+\int\limits_{[0,\infty)}T(I+sT)xd\lambda(s)\quad (x\in D(T)).
$$
Since
$$
f(w)=a+bw+\int\limits_{(0,\infty)}\frac{w}{s^{-1}+w}s^{-1}d\lambda(s),
$$
every such function is complete Bernstein. Consequently, the functional calculus under consideration is in fact a form of Hirsch functional calculus.

In what follows ($\mathcal{I}, \|\cdot\|_{\mathcal{I}}$) stands for a \textit{symmetrically normed operator ideal} in complex Banach space $X,$ i.e., two-sided ideal of the algebra $\mathcal{L}(X)$
of bounded operators on  $X$ that is complete with respect to the norm $\|\cdot\|_{\mathcal{I}},$ and satisfies $\|ASB\|_{\mathcal{I}}\leq \|A\||S\|_{\mathcal{I}}\|B\|,$  and $\|S\|\leq \||S\|_{\mathcal{I}}$
for all $A,B\in \mathcal{L}(X)$ and $S\in  \mathcal{I}$ (the case $\mathcal{I} = \mathcal{L}(X)$ is also of interest).

\textbf{3. Operator differentiability}

\textbf{Definition 4}.  Let $\varphi\in\mathcal{ OM}_-,$ $A \in \mathcal{NP}(X),$ and let $\mathcal{I}$ be an operator ideal. A bounded operator $\varphi^\nabla_A$
on $\mathcal{I}$ (transformer) is called an $\mathcal{I}$\textit{-Frechet derivative} of a function $\varphi$ at the point $A$ if for every $\Delta A\in \mathcal{I}$ the following
asymptotic equality holds:
$$
\|\varphi(A+\Delta A)-\varphi(A)-\varphi_A^\nabla(\Delta A)\|_\mathcal{I}=o(\|\Delta A\|_\mathcal{I})
 \mbox{ as }  \|\Delta A\|_\mathcal{I}\to 0.
$$

We need the following  lemmas.

\textbf{Lemma 1.} (i) \textit{Let $A \in \mathcal{N}(X).$ Then  $A+V \in \mathcal{N}(X)$ and $\rho(A+V)\supset \rho(A)$ for every $V\in \mathcal{L}(X)$ such that $\|V\|<1/M_A.$ In this case $M_{A+V}\leq M_A/(1-M_A\|V\|).$}

(ii) \textit{Let $A \in \mathcal{NP}(X).$ Then  $A+V \in \mathcal{NP}(X)$ and $\rho(A+V)\supset \rho(A)$ for every $V\in \mathcal{L}(X)$ such that $\|V\|<1/2M_A.$ In this case  $M_{A+V}\leq 2M_A/(1-2M_A\|V\|).$}

Proof. (i) First note that   $\|V\|<1/M_A\leq 1/\|R(t,A)\|$ for all $t\in \mathbb{R}_+.$ It follows in view of \cite[Remark IV.3.2]{Kato} that $\rho(A+V)\supset \rho(A)\supset\mathbb{ R}_+.$

Next, applying \cite[Theorem IV.1.16, Remark IV.1.17]{Kato} we have for  $t\in \mathbb{R}_+$
$$
\|R(t,A+V)-R(t,A)\|\leq \frac{\|V\|\|R(t,A)\|^2}{1-\|V\|\|R(t,A)\|}\eqno(3)
$$
since $\|V\|\|R(t,A)\|< \|R(t,A)\|/M_A\leq 1.$

Thus,   we obtain for  $t\in \mathbb{R}_+$
$$
\|R(t,A+V)\|\leq \|R(t,A)\|+\frac{\|V\|\|R(t,A)\|^2}{1-\|V\|\|R(t,A)\|}=\frac{\|R(t,A)\|}{1-\|V\|\|R(t,A)\|}
$$
$$
\leq
\frac{\|R(t,A)\|}{1-\|V\|M_A}\leq \frac{M_A/(1-M_A\|V\|)}{1+t}
$$
which completes the proof of (i).

(ii) For all  $\varepsilon>0$ we have $M_{A-\varepsilon I}\leq 2M_A,$  since
$$
M_{A-\varepsilon I}=\sup\limits_{\lambda>0}\|\lambda((\lambda+\varepsilon)I-A)^{-1}\|\leq
 \sup\limits_{\lambda>0}\|(\lambda+\varepsilon)((\lambda+\varepsilon)I-A)^{-1}\|+
 $$
 $$
 \varepsilon\sup\limits_{\lambda>0}\|((\lambda+\varepsilon)I-A)^{-1}\|\leq M_A+\varepsilon\sup\limits_{\lambda>0}\frac{M_A}{\lambda+\varepsilon}=2M_A.
$$
Therefore for all  $\varepsilon>0$ the condition  $\|V\|<1/2M_A$ yields  $\|V\|<1/M_{A-\varepsilon I}.$ So, by (i), $A-\varepsilon I +V\in \mathcal{N}(X).$ It follows first of all that  $\rho(A-\varepsilon I +V)\supset \rho(A-\varepsilon I)$ which implies $\rho(A +V)\supset \rho(A).$

Moreover,
$$
\sup_{t>0}\|tR(t,(A-\varepsilon I) +V)\|\leq M_{A-\varepsilon I}/(1-M_{A-\varepsilon I}\|V\|)\leq 2M_A/(1-2M_A\|V\|).
$$
Consequently, for all $\varepsilon>0, t>0$
$$
\|tR(t+\varepsilon, A+V)\|\leq 2M_A/(1-2M_A\|V\|).
$$
Letting $\varepsilon$   tend to zero, we obtain
$\|tR(t, A+V)\|\leq 2M_A/(1-2M_A\|V\|)$ for all $t>0$ and the result follows.

\textbf{Lemma 2.} (Cf. \cite{Nab}.) (i) \textit{Let $\varphi\in \mathcal{OM}_-$. For any operators $A,B\in \mathcal{N}(X)$ such that $D(A)\subseteq D(B),$ and  $A- B\in \mathcal{I}$ the operator $\varphi(A)-\varphi(B)$ belongs to $\mathcal{I},$ too, and satisfies the inequality}
$$
\|\varphi(A)-\varphi(B)\|_{\mathcal{I}}\leq M_AM_B\varphi'(-1)\|A-B\|_{\mathcal{I}}.
$$

(ii) \textit{If, in addition, $\varphi'(-0)\ne\infty,$ then $\varphi(A)-\varphi(B)$ belongs to $\mathcal{I}$ for any operators $A,B\in \mathcal{NP}(X)$ such that $D(A)\subseteq D(B),$ and  $A- B\in \mathcal{I}$ and}
$$
\|\varphi(A)-\varphi(B)\|_{\mathcal{I}}\leq M_AM_B\varphi'(-0)\|A-B\|_{\mathcal{I}}.
$$

Proof. (i) Since $AR(t,A)x=R(t,A)Ax$ for $x\in D(A),$ we have
$$
 (\varphi(A)-\varphi(B))x=\int\limits_{(0,\infty)}(AR(t,A)-BR(t,B))xd\mu(t)\ (x\in D(A))
$$
Let $g(t):=AR(t,A)-BR(t,B)\ (t>0).$  The well known equality
$AR(t,A)=-I+tR(t,A)\quad (t\in\rho(A))$
implies in view of the second resolvent identity that
$g(t)=t(R(t,A)-R(t,B))=tR(t,A)(A-B)R(t,B).$
Therefore  $\|g(t)\|_{\mathcal{I}}\leq M_AM_B\|A-B\|_{\mathcal{I}}t/(1+t)^2.$
It follows that the Bochner integral $\int_{(0,\infty)}g(t)d\mu(t)$ exists with respect to the $\mathcal{I}$ norm, and the desired inequality  is valid, because $\int_{(0,\infty)}td\mu(t)/(1+t)^2=\varphi'(-1).$

(ii) In this case, $\|g(t)\|_{\mathcal{I}}\leq M_AM_B\|A-B\|_{\mathcal{I}}/t.$ Then the Bochner integral $\int_{(0,\infty)}g(t)d\mu(t)$ exists with respect to the $\mathcal{I}$ norm and the desired inequality holds, since  $\int_{(0,\infty)}d\mu(t)/t=\varphi'(-0)\ne\infty.$

\textbf{Definition 5}. We introduce on $\mathcal{NP}(X)$ the following equivalence relation: operators $A$ and $A'$ from  $\mathcal{NP}(X)$ are
equivalent if $A'-A\in \mathcal{L}(X)$. Formula $\|A'-A\|$  defines metrics in every equivalence class.

\textbf{Theorem 1.} 1) (a) \textit{A function $\varphi\in\mathcal{ OM}_-$ is $\mathcal{I}$-differentiable in Frechet sense at any point $A\in \mathcal{N}(X),$ and its $\mathcal{I}$-Frechet derivative  is given by the formula}
$$
\varphi_A^\nabla(V)=\int\limits_0^\infty R(t,A)VR(t,A)td\mu(t)\quad (V\in \mathcal{I}). \eqno(4)
$$

(b)  \textit{For every equivalence class $\mathcal{C}$ of operators from  $\mathcal{N}(X)$ the mapping $A\mapsto  \varphi_A^\nabla$ from $\mathcal{C}$ to $\mathcal{L}(\mathcal{L}(X))$ is continuous.}

2) \textit{Let the function $\varphi\in\mathcal{ OM}_-$ satisfies $\varphi'(-0)\ne\infty, \varphi''(-0)\ne\infty.$ Then $\varphi$ is $\mathcal{I}$-differentiable in Frechet sense at any point $A\in \mathcal{NP}(X),$ and its $\mathcal{I}$-Frechet derivative  is given by formula (4).}

Proof. 1) (a). Fix $A\in \mathcal{N}(X).$ The transformer $F_A:\mathcal{I}\to \mathcal{I}$ which is defined by the right-hand side of formula (4) is bounded, because
$$
\|F_A(V)\|_{\mathcal{I}}\leq M_A^2 \int\limits_0^\infty \frac{td\mu(t)}{(1+t)^2}\|V\|_{\mathcal{I}}.
$$

Let $\Delta A \in \mathcal{I}$ be such that $\|\Delta A\|_{\mathcal{I}}<1/M_A.$  Then $A+\Delta A \in \mathcal{N}(X)$  (see lemma 1) and
by means of consideration from the proof of lemma 2
 we obtain the equality
$$
\varphi(A+\Delta A)-\varphi(A)-F_A(\Delta A)=\int\limits_0^\infty(R(t,A+\Delta A)-R(t,A))\Delta AR(t,A)td\mu(t).\eqno(5)
$$
Furthermore,  we have by lemma 1
$$
\|R(t,A+\Delta A)-R(t,A)\|=\|R(t,A+\Delta A)\Delta AR(t,A)\|\leq
$$
$$
\frac{M_{A+\Delta A}M_A}{(1+t)^2}\|\Delta A\|_{\mathcal{I}}\leq \frac{2M_A^2}{1-M_A\|\Delta A\|_{\mathcal{I}}}\frac{\|\Delta A\|_{\mathcal{I}}}{(1+t)^2}.
$$
In view of  this inequality formula (5) implies
$$
\|\varphi(A+\Delta A)-\varphi(A)-F_A(\Delta A)\|_{\mathcal{I}}\leq 2M_A^3\int\limits_0^\infty\frac{td\mu(t)}{(1+t)^3}\frac{\|\Delta A\|_{\mathcal{I}}^2}{1-2M_A\|\Delta A\|_{\mathcal{I}}}=o(\|\Delta A\|_{\mathcal{I}})
$$
and the first statement follows.

(b). Let operators $A$ and $A'$ from  $\mathcal{N}(X)$ be equivalent. By virtue of formula (4) for any $B\in \mathcal{L}(X)$ we have
$$
(\varphi^\nabla_{A'}-\varphi^\nabla_A)(B)=\int\limits_0^\infty R(t,A')B(R(t,A')-R(t,A))td\mu(t)+
\int\limits_0^\infty(R(t,A')-R(t,A))BR(t,A)td\mu(t).
$$
Hence,
$$
\|(\varphi^\nabla_{A'}-\varphi^\nabla_A)(B)\|\leq 2\max\{M_{A'},M_A\}\|B\|\int\limits_0^\infty\|R(t,A')-R(t,A)\|\frac{t}{1+t}d\mu(t).\eqno(6)
$$
Choose arbitrary $\varepsilon\in (0,1)$ and let $\|A'-A\|<\varepsilon/(2M_A).$ Then $\|A'-A\|\|R(t,A)\|<1/2,$ and therefore  $M_{A'}\leq 2M_A$ by lemma 1 with $V=A'-A.$ Moreover,  formula (3) implies for $V=A'-A$ that
$$
\|R(t,A')-R(t,A)\|<2\|A'-A\|\|R(t,A)\|^2\leq 2\|A'-A\|\frac{M_A^2}{(1+t)^2}<\frac{M_A}{(1+t)^2}\varepsilon.
$$
So, by virtue of formula (6)
$$
\|(\varphi^\nabla_{A'}-\varphi^\nabla_A)B\|\leq 4M_A^3\int\limits_0^\infty\frac{td\mu(t)}{(1+t)^3}\varepsilon \|B\|,
$$
and then
$$
\|\varphi^\nabla_{A'}-\varphi^\nabla_A\|_{\mathcal{L}(\mathcal{L}(X))}\leq \left( 4M_A^3\int\limits_0^\infty\frac{td\mu(t)}{(1+t)^3}\right)\varepsilon ,
$$
which completes the proof of the part (b).

2). Since $\int_{(0,\infty)}d\mu(t)/t=\varphi'(-0)\ne\infty,$ the transformer $F_A:\mathcal{I}\to \mathcal{I}$ which is defined by the right-hand side of formula (4) is bounded. Moreover, formula (5) implies by the second resolvent identity that
$$
\|\varphi(A+\Delta A)-\varphi(A)-F_A(\Delta A)\|_{\mathcal{I}}=
$$
$$
\left\|\int\limits_0^\infty R(t,A+\Delta A)\Delta AR(t,A)\Delta AR(t,A)td\mu(t)\right\|_{\mathcal{I}}\leq \int\limits_0^\infty\|R(t,A+\Delta A)\|\|R(t,A)\|^2td\mu(t)\|\Delta A\|_{\mathcal{I}}^2.
$$
If $\|\Delta A\|_{\mathcal{I}}<1/2M_A$ we have by lemma 1 that $\|R(t,A+\Delta A)\|\leq M_{A+\Delta A}/t\leq 2M_A(1-2M_A\|\Delta A\|)^{-1}/t.$ Since $\int_{(0,\infty)}d\mu(t)/t^2=\varphi''(-0)\ne\infty,$ it follows that
$$
\|\varphi(A+\Delta A)-\varphi(A)-F_A(\Delta A)\|_{\mathcal{I}}\leq 2M_A^3\int\limits_0^\infty\frac{d\mu(t)}{t^2}\frac{\|\Delta A\|_{\mathcal{I}}^2}{1-2M_A\|\Delta A\|_{\mathcal{I}}}=o(\|\Delta A\|_{\mathcal{I}})
$$
and the proof is complete.

\textbf{Remark 4}. If $\varphi\in\mathcal{ OM}_-$ has integral representation (1) (with $a=b=0$), then $\varphi'(s)=\int_0^\infty td\mu(t)/(t-s)^2.$ Therefore for $A\in \mathcal{N}(X)$ or $A\in \mathcal{NP}(X)$ and $\varphi'(-0)\ne \infty$ one can put
$$
\varphi'(A):=\int\limits_0^\infty R(t,A)^2td\mu(t).
$$
If $V\in\mathcal{I}$ and $A$ and $V$ commutes it is easy to verify that $\varphi^\nabla_A(V)=\varphi'(A)V.$

\textbf{Examples.} 1) The function
$\psi_\lambda(s):=\log \lambda-\log(\lambda-s)$ ($\lambda>0$) belongs to $\mathcal{OM}_-$ and has the integral representation ($s<0$)
$$
\psi_\lambda(s)=\int\limits_\lambda^\infty\frac{s}{t-s}\frac{dt}{t}.
$$
Thus, theorem 1 implies that for $A\in \mathcal{NP}(X),$ and $V\in\mathcal{I}$
$$
(\psi_\lambda)^\nabla_A(V)=\int\limits_\lambda^\infty R(t,A)VR(t,A)dt.
$$
If $A$ and $V$ commutes it follows that $(\psi_\lambda)^\nabla_A(V)=\psi_\lambda'(A)V=(\lambda I-A)^{-1}V.$

2) The function
$\varphi(s):=-(-s)^\alpha$ ($\alpha\in(0,1)$) belongs to $\mathcal{OM}_-$ and has the integral representation
$$
\varphi(s)=\frac{\sin\alpha\pi}{\pi}\int\limits_0^\infty\frac{s}{t-s}t^{\alpha-1}dt.
$$
So, theorem 1 implies that for $A\in \mathcal{N}(X),$ and $V\in\mathcal{I}$

$$
\varphi^\nabla_A(V)=\frac{\sin\alpha\pi}{\pi}\int\limits_0^\infty R(t,A)VR(t,A)t^\alpha dt.
$$

Before  we formulate our   next theorem note that there exist differentiable functions $f$ on $\mathbb{R}$ such that for some self-adjoint operators $A$ and $V$  the function $t\mapsto f(A + tV)-f(A)$ is not differentiable at the origin (see, e. g., \cite[Theorem 1.2.8]{AP1}).

\textbf{Theorem 2.} \textit{If $\varphi\in\mathcal{ OM}_-,$ $A \in \mathcal{N}(X),$ $V\in \mathcal{I}$, then the $\mathcal{I}$-valued function $z\mapsto \varphi(A + zV)-\varphi(A)$
is analytic in the neighborhood of the origin $\mathcal{O}_{A,V} := \{z \in \mathbb{C} : |z| <1/(\|V\|_{\mathcal{I}}M_A) \}$, and it allows in this
neighborhood the expansion}
$$
\varphi(A+zV)-\varphi(A)=\sum\limits_{n=1}^\infty z^n C_n,\eqno(7)
$$
\textit{where the series absolutely converges in the $\mathcal{I}$ norm and}
 $$
 C_n=\frac{1}{n!}\left.\frac{d^n}{dz^n}\varphi(A+zV)\right|_{z=0}=\int\limits_0^\infty(R(t,A)V)^nR(t,A)td\mu(t)\eqno(8)
 $$
(\textit{the derivatives are understood in the sense of the $\mathcal{I}$ norm}).

Proof. For every $z\in \mathcal{O}_{A,V}$  operator $A+zV$ belongs to $\mathcal{N}(X)$ by lemma 1 and therefore  $\varphi(A+zV)-\varphi(A)\in \mathcal{I}$ by lemma 2. Since the function $\varphi$ is $\mathcal{I}$-Frechet differentiable at the point $A+zV$ by theorem 1, its   $\mathcal{I}$-Gateaux derivative  at the point $A+zV,$ the transformer $d/dh\varphi(A+(z+h)V)|_{h=0},$    coincides with  $\varphi^\nabla_{A+zV}.$ This means, due to formula (4), that
$$
\frac{d}{dz}\varphi(A+zV)=\varphi^\nabla_{A+zV}=\int\limits_0^\infty R(t,A+zV)VR(t,A+zV)td\mu(t). \eqno(9)
$$
Consequently, $\mathcal{I}$-valued function $\varphi(A + zV)-\varphi(A)$ is analytic in  $\mathcal{O}_{A,V}$ and allows an
expansion (7), where $C_n$ is determined by the first of equalities (8).The second equality is the  consequence of the following  equality $(z\in \mathcal{O}_{A,V})$
$$
\frac{d^n}{dz^n}\varphi(A+zV)=n!\int\limits_0^\infty(R(t,A+zV)V)^nR(t,A+zV)td\mu(t)\eqno(10)
 $$
which we will prove  by induction. For $n = 1$ it holds by virtue of (9). Assume that it is valid for certain $n$ and let $|z|<q/(\|V\|_{\mathcal{I}}M_A)$ for $q\in (0,1)$. Since $d/dzR(t,A+zV)=R(t,A+zV)VR(t,A+zV),$ we have differentiating under the integral sign
$$
\frac{d^{n+1}}{dz^{n+1}}\varphi(A+zV)=n!\int\limits_0^\infty\frac{d}{dz}((R(t,A+zV)V)^nR(t,A+zV))td\mu(t)=
$$
$$
n!\int\limits_0^\infty((nR(t,A+zV)V)^{n-1}\frac{d}{dz}(R(t,A+zV))VR(t,A+zV)+(R(t,A+zV)V)^n\frac{d}{dz}R(t,A+zV))td\mu(t)= $$
$$
(n+1)!\int\limits_0^\infty(R(t,A+zV)V)^{n+1}R(t,A+zV)td\mu(t).
 $$
Since,   by lemma 1, $M_{A+zV}\leq M_A/(1-M_A\|zV\|_{\mathcal{I}})<M_A/(1-q),$ the validity of differentiation under the integral sign follows from the estimate
$$
\|(R(t,A+zV)V)^{n+1}R(t,A+zV)t\|_{\mathcal{I}} \leq  \|R(t,A+zV)\|^{n+2}\|V\|_{\mathcal{I}}^{n+1}t\leq
$$
$$
\left(\frac{M_{A+zV}}{1+t}\right)^{n+2}\|V\|_{\mathcal{I}}^{n+1}t<
\left(\frac{M_A}{1-q}\right)^{n+2}\|V\|_{\mathcal{I}}^{n+1}\frac{t}{(1+t)^{n+2}}.
$$
Finally, for $z\in \mathcal{O}_{A,V}$ we have
$$
\|z^nC_n\|_{\mathcal{I}}\leq \frac{1}{\|V\|_{\mathcal{I}}^{n}M_A^n}\|C_n\|_{\mathcal{I}}\leq
\frac{1}{\|V\|_{\mathcal{I}}^{n}M_A^n}\int\limits_0^\infty\|(R(t,A)V)^nR(t,A)\|_{\mathcal{I}}td\mu(t)\leq
M_A\int\limits_0^\infty\frac{td\mu(t)}{(1+t)^{n+1}}.
$$
Since
$$
\sum\limits_{n=1}^\infty\int\limits_0^\infty\frac{td\mu(t)}{(1+t)^{n+1}}=\int\limits_0^\infty\frac{td\mu(t)}{(1+t)^{2}}<\infty,
$$
 the series in (7) absolutely converges in the $\mathcal{I}$ norm. This completes the proof.

The analog of theorem 2 is valid for $A\in \mathcal{NP}(X)$ as well.

\textbf{Theorem 3.} \textit{Let $\varphi\in\mathcal{ OM}_-$ be such that $\varphi^{(n)}(-0)\ne\infty$ for all $n\in \mathbb{N},$ $A \in \mathcal{NP}(X),$ $V\in \mathcal{I}.$ Then the $\mathcal{I}$-valued function $z\mapsto \varphi(A + zV)-\varphi(A)$
is analytic in the neighborhood of the origin $\mathcal{O}_{A,V}^{(2)} := \{z \in \mathbb{C} : |z| <1/(2\|V\|_{\mathcal{I}}M_A) \}$, and it allows in this
neighborhood the expansion (7)
where the series absolutely converges in the $\mathcal{I}$ norm and formula (8) holds.}

 The proof is similar to the proof of theorem 2.

\textbf{4. Application to   perturbation determinants}

We shall apply theorem 3 to prove the important formula for a logarithmic derivative of a perturbation determinant   of  operators on   Banach spaces  \cite{OaMII}, \cite{arX2018}. Recall that operator on $X$ is \textit{ nuclear} if it is representable as the sum of absolutely convergent in operator norm series  of rank one operators; if, in addition, $X$ has the approximation property, the continuous trace $\mathrm{tr}$ is defined on the operator ideal $\mathfrak{S}_1=\mathfrak{S}_1(X)$ of nuclear operators on $X,$ see, e.g., \cite{DF}.

\textbf{Definition 6}   \cite{arX2018}. Let $X$ has the approximation property.  For $A, B\in \mathcal{NP}(X)$ such that $D(A)\subseteq D(B)$ and $B-A$  is nuclear define \textit{the perturbation determinant  for the  pair} $(A,B)$  by
$$
\Delta_{B/A}(\lambda):=\exp(\mathrm{tr}(\psi_\lambda(B)-\psi_\lambda(A)))\quad (\lambda>0)
$$
where $\psi_\lambda(s):=\log \lambda-\log(\lambda-s)$ (see example 1 above).

This is  a generalization of the classical notion of a perturbation determinant  for  pairs of self-adjoint operators on Hilbert space (see, e.g., \cite[formula (3.25)]{BY}).

Now we list several important properties of a perturbation determinant.

1) If operators $A, B, C\in \mathcal{NP}(X)$ be such that $A-B$ and $B-C$ are nuclear, then
$$
\Delta_{B/A}(\lambda)\Delta_{C/B}(\lambda)=\Delta_{C/A}(\lambda),
$$

in particular,
$$
\Delta_{B/A}(\lambda)\Delta_{A/B}(\lambda)=1\quad (\lambda>0).
$$
This property  is an immediate consequence of the definition.

2)  $\lim\limits_{\lambda\to+\infty}\Delta_{B/A}(\lambda)=1.$

Indeed, by the lemma 2
$$
\|\psi_\lambda(A)-\psi_\lambda(B)\|_{\mathfrak{S}_1}\leq M_AM_B\psi'_\lambda(-0)\|A-B\|_{\mathfrak{S}_1}=\frac{1}{\lambda}M_AM_B\|A-B\|_{\mathfrak{S}_1}\to 0 \mbox{ as } \lambda\to+\infty.
$$

 Since  $\psi_\lambda$ enjoys the conditions of  theorem 3, we get also the next property.

\textbf{Corollary 1.} \textit{Let $X$ has the approximation property.  Let $A\in \mathcal{NP}(X)$ and $V\in \mathfrak{S}_1.$  Then the map $z\mapsto \Delta_{(A+zV)/A}(\lambda)$   is analytic in $\mathcal{O}_{A,V}^{(2)}$ for every $\lambda>0$ and it allows in this
neighborhood the representation}
$$
\Delta_{(A+zV)/A}(\lambda)=\prod\limits_{n=1}^\infty e^{\mathrm{tr}(C_n(\lambda)) z^n}\quad (\lambda>0)
$$
\textit{where }
$$
\mathrm{tr}(C_n(\lambda))=\int\limits_\lambda^\infty\mathrm{tr}((R(t,A)V)^nR(t,A))dt.
$$

Proof. This follows from theorem 3 applied to the function $\psi_\lambda.$

 Now we are able to prove  a generalization of the   formula for the logarithmic derivative of a perturbation determinant  (for negative operators on Banach spaces this formula was proved by another  method in \cite{arX2018}).

\textbf{Theorem 4.} \textit{Let $X$ has the approximation property, $A,B\in \mathcal{NP}(X), A\ne B$ and $B-A\in \mathfrak{S}_1(X).$ Then for all $\lambda>1/2\|A-B\|$}
$$
\frac{\Delta'_{B/A}(\lambda)}{\Delta_{B/A}(\lambda)}=\mathrm{tr}(R(\lambda,A)-R(\lambda,B)).
$$

Proof. The proof of this theorem is broken into two steps.

1) First consider the case $\|B-A\|_{\mathfrak{S}_1}<1/2M_A.$ One can assume that $B-A=zV,$ where $V\in \mathfrak{S}_1(X),$ $z\in \mathcal{O}_{A,V}^{(2)}$ (for example, we can take $V=(B-A)/\|B-A\|_{\mathfrak{S}_1},$ $z=\|B-A\|_{\mathfrak{S}_1}$). Then, by corollary 1,
$$
\frac{\Delta'_{B/A}(\lambda)}{\Delta_{B/A}(\lambda)}=\frac{d}{d\lambda}\log\Delta_{(A+zV)/A}(\lambda)=
\sum\limits_{n=1}^\infty\frac{d}{d\lambda}\mathrm{tr}(C_n(\lambda)) z^n\  (\lambda>M_A)  \eqno(11)
$$
because the series in the right-hand side converges uniformly in $\lambda>M_A.$ Indeed,
by corollary 1
$$
\frac{d}{d\lambda}\mathrm{tr}(C_n(\lambda))=-\mathrm{tr}((R(\lambda,A)V)^nR(\lambda,A))\eqno(12)
$$
and for all $\lambda>M_A$ and $z\in \mathcal{O}_{A,V}^{(2)}$ we have
$$
 \|zR(\lambda,A)V\|\leq \|zR(\lambda,A)V\|_{\mathfrak{S}_1}\leq |z|\frac{M_A}{\lambda}\|V\|_{\mathfrak{S}_1}<\frac{1}{2}. \eqno(13)
$$
Consequently,
$$
\left|\frac{d}{d\lambda}\mathrm{tr}(C_n(\lambda))z^n\right|=|\mathrm{tr}((zR(\lambda,A)V)^nR(\lambda,A))|\leq \|(zR(\lambda,A)V)^nR(\lambda,A)\|_{\mathfrak{S}_1}=
$$
$$
\|(zR(\lambda,A)V)^{n-1}(zR(\lambda,A)V)R(\lambda,A)\|_{\mathfrak{S}_1}\leq \|(zR(\lambda,A)V)^{n-1}\|\|zR(\lambda,A)V\|_{\mathfrak{S}_1}\|R(\lambda,A)\|\leq
$$
$$
\frac{1}{2^{n-1}}|z|\|V\|_{\mathfrak{S}_1}\left(\frac{M_A}{\lambda}\right)^2 <\frac{1}{M_A}\frac{1}{2^n}.\eqno(14)
$$
Formulas (11) and (12) imply that
$$
\frac{\Delta'_{B/A}(\lambda)}{\Delta_{B/A}(\lambda)}=-\mathrm{tr}
\sum\limits_{n=1}^\infty ((zR(\lambda,A)V)^nR(\lambda,A))\quad  (\lambda>M_A)  \eqno(15)
$$
(the series in the right-hand side converges absolutely with respect to the nuclear norm as the proof of (14) shows).

Moreover, since the estimate (13) is valid for $\lambda>M_A$ and $z\in \mathcal{O}_{A,V}^{(2)},$ we have
$$
\sum\limits_{n=1}^\infty ((zR(\lambda,A)V)^n=zR(\lambda,A)V(I-zR(\lambda,A)V)^{-1}.  \eqno(16)
$$

We claim that
$$
zR(\lambda,A)V(I-zR(\lambda,A)V)^{-1}R(\lambda,A)=R(\lambda,A+zV)-R(\lambda, A). \eqno(17)
$$
To this end, we shall show that for all  $x\in D(A)$
$$
(I-zR(\lambda,A)V)^{-1}x=R(\lambda,A+zV)(\lambda I-A)x. \eqno(18)
$$
Indeed, this follows from the next calculations ($\mathrm{Im}(R(\lambda,A+zV))=D(A)$):
$$
(I-zR(\lambda,A)V)R(\lambda,A+zV)(\lambda I-A)x=(R(\lambda,A)(\lambda I-A)-zR(\lambda,A)V)R(\lambda,A+zV)(\lambda I-A)x=
$$
$$
R(\lambda,A)(\lambda I-A-zV)R(\lambda,A+zV)(\lambda I-A)x=R(\lambda,A)(\lambda I-A)x=x.
$$
Now, in view  of (18) the left-hand side of (17) takes the form
$$
zR(\lambda,A)V(I-zR(\lambda,A)V)^{-1}R(\lambda,A)=zR(\lambda,A)VR(\lambda,A+zV)(\lambda I-A)R(\lambda,A)=
$$
$$
zR(\lambda,A)VR(\lambda,A+zV)  =-(R(\lambda, A)-R(\lambda,A+zV))
 $$
and (17) follows.
Putting  together (17), (16) and (15) we get for  $\lambda>M_A$ that
$$
\frac{\Delta'_{B/A}(\lambda)}{\Delta_{B/A}(\lambda)}=\mathrm{tr}(R(\lambda, A)-R(\lambda,A+zV))=
\mathrm{tr}(R(\lambda,A)-R(\lambda,B)).
$$
For $\lambda>M_A$ the proof of the case 1)  is complete.

2) Now let  $\|B-A\|_{\mathfrak{S}_1}$ be an arbitrary positive number, and $\lambda>1/2\|A-B\|.$ Note that for all positive $k$ and $\lambda$ we have $R(\lambda, kA)=k^{-1}R(k^{-1}\lambda, A).$ First of all it follows that $kA, kB\in \mathcal{NP}(X)$ and $M_{kA}=M_A.$ Now we choose $k<1/(2M_A\|A-B\|$) (and then $\|kA-kB\|<1/2M_A$). According to the case 1) we have for $\lambda>M_A$
$$
\frac{d}{d\lambda}(\log\Delta_{kB/kA})(\lambda)=
\mathrm{tr}(R(\lambda,kA)-R(\lambda,kB))=k^{-1}\mathrm{tr}(R(k^{-1}\lambda,A)-R(k^{-1}\lambda,B)).
$$
On the other hand,
$$
\log\Delta_{kB/kA}(\lambda)=\mathrm{tr}(\psi_\lambda(kB)-\psi_\lambda(kA))=\mathrm{tr}(\psi_{k^{-1}\lambda}(B)-\psi_{k^{-1}\lambda}(A))=
\log\Delta_{B/A}(k^{-1}\lambda)
$$
which implies
$$
\frac{d}{d\lambda}(\log\Delta_{kB/kA})(\lambda)=k^{-1}\frac{d}{d\lambda}(\log\Delta_{B/A})(k^{-1}\lambda).
$$
Thus the equality under consideration is valid for all $\lambda>kM_A,$ the more so, this is true for $\lambda>1/2\|A-B\|.$ This completes the proof.

\textbf{Remark.} Using the same arguments as in the proof of theorem 4.2 in \cite{arX2018} (see formula (10) there) one can show that $\Delta_{A,B}$ possesses an analytic continuation into some open sector $S_\theta\subset \rho(A)\cap \rho(B)$ symmetric about the positive real semiaxis. It follows that the above mentioned properties of the perturbation determinant are valid for $\lambda\in S_\theta.$

\textbf{Corollary 2.}  \textit{Let the conditions of theorem 4 are fulfilled. Let the Banach space $X$ has the extra property that the trace   on  $\mathfrak{S}_1(X)$ is nilpotent in a sense  that  $\mathrm{tr}(N)=0$ for every nilpotent operator $N.$ Suppose $z_1$ is a regular point or an isolated eigenvalue of the operators
$B$ and $A$ of finite algebraic multiplicities $k_0$ and $k.$ Then at the point $z_1$ the
function $\Delta_{B/A}(z)$ has a pole (or zero) of order $k_0 - k$ (respectively of order $k - k_0$).}

Proof. Due to the theorem 4 the proof of this assertion is similar to the proof of the property 4 of the perturbation determinant given in \cite[p. 267]{Ya} (formula (5) from \cite[pp. 266--267]{Ya} is valid for Banach spaces,  see \cite[Chapter III, subsection 6.5]{Kato}).

\textbf{Acknowledgments.
}
This work was financially supported by the  Fund of Fundamental Research of Republic of Belarus. Grant number $\Phi$ 17-082.

\begin{flushleft}
\textsc{Department of mathematics and programming technologies}\\

\textsc{ F. Skorina Gomel State University}\\

\textsc{104 Sovietskaya St.,
 Gomel, 246019, Belarus}\\

\textsc{E-mail}: amirotin@yandex.ru
\end{flushleft}
\end{document}